\documentclass[titlepage,11pt]{article}
\usepackage{latexsym}
\usepackage{amsfonts}
\usepackage{amsmath}
\newcommand\blackslug{\hbox{\hskip 1pt \vrule width 4pt height 8pt depth 1.5pt
        \hskip 1pt}}
\newcommand\bbox{\hfill \quad \blackslug \bigbreak}
\def\d{\hbox{-}}
\def\c{\hbox{-}\cdots\hbox{-}}
\def\l{,\ldots,}

\title{Induced subgraphs of graphs with large chromatic number. 
\\I. Odd holes}
\author{Alex Scott\\
Mathematical Institute, University of Oxford, Oxford OX2 6GG, UK
\\
\\
Paul Seymour\thanks{Supported by ONR grant N00014-10-1-0680 and NSF
grant DMS-1265563.}\\
Princeton University, Princeton, NJ 08544, USA}

\date{August 19, 2014; revised \today}

\newtheorem{thm}{}[section]

\newcommand{\Proof}{\noindent{\bf Proof.}\ \ }

\begin{document}
\maketitle
\begin{abstract}
An {\em odd hole} in a graph is an induced subgraph which is a cycle of odd length at least five.
In 1985, A. Gy\'arf\'as made the conjecture that for all $t$ there exists $n$ such that every graph with no $K_t$ subgraph and no odd hole is $n$-colourable.
We prove this conjecture.
\end{abstract}

\section{Introduction}

All graphs in this paper are finite and have no loops or parallel edges. A {\em hole} in a graph is an induced subgraph which is a cycle of length
at least four, and a hole is {\em odd} if its length is odd.  An {\em odd antihole} is an induced subgraph that is the complement of an odd hole.
We denote the chromatic number of a graph $G$ by $\chi(G)$. A {\em clique} of $G$ is a subset of $V(G)$ such that all its members are pairwise adjacent, 
and the cardinality of the largest clique of $G$ is denoted
by $\omega(G)$. A clique of cardinality $W$ is called a {\em $W$-clique}. 

Clearly $\chi(G)\ge \omega(G)$, but $\chi(G)$ cannot in general be bounded above by a function of $\omega(G)$; indeed,
there are graphs with $\omega(G)=2$ and $\chi(G)$ arbitrarily large~\cite{tutte1,tutte2}. 
On the other hand,
the ``strong perfect graph theorem''~\cite{CRST} says: 

\begin{thm}\label{pgt} 
Let $G$ be a graph that does not contain an odd hole and does not contain an odd antihole.  Then
$\chi(G)=\omega(G)$.
\end{thm}

This is best possible in a sense, since odd holes
and odd antiholes themselves have $\chi>\omega$.
But what happens if we only exclude odd holes? A. Gy\'arf\'as~\cite{gyarfas} proposed the following
conjecture in 1985, and it is the main result of this paper:

\begin{thm}\label{mainconj}
{\bf (Conjecture)} There is a function $f$ such that
$\chi(G)\le f(\omega(G))$ for every graph $G$ with no odd hole.
\end{thm}

Let $G$ have no odd hole. 
It is trivial that $\chi(G)=\omega(G)$ if $\omega(G)\le 2$, and it was proved 
in~\cite{4oddhole} that $\chi(G)\le 4$ if $\omega(G)=3$;
but for larger values of $\omega(G)$, no bound on $\chi(G)$ in terms of $\omega(G)$ was known.
Here we prove conjecture \ref{mainconj}, and give an explicit function $f$.

\begin{thm}\label{mainthm}
Let $G$ be a graph with no odd hole.  Then $\chi(G)\le 2^{2^{\omega(G)+2}}$.
\end{thm}

We imagine that this bound is nowhere near best possible. There is an example that
satisfies
$\chi(G) \ge
\omega(G)^\alpha$, where $\alpha=\log(7/2)/\log(3)\approx 1.14$, but this is the best we have found.
To see this, let $G_0$ have one vertex, and for $k\ge 1$ let $G_k$ be obtained
from $G_{k-1}$ by substituting a seven-vertex antihole for each vertex. Then 
$G_k$ has no odd hole, $\omega(G_k)=3^k$, and 
$\chi(G_k)\ge (7/2)^k$ (because the largest stable set in $G_k$ has 
cardinality $2^k$).

\section{Background}

What can we say about the structure of a graph $G$ with large chromatic number, if it does not contain a large clique?
In particular, what can we say about its induced subgraphs?
If $\chi(G)>\omega(G)$, then the strong perfect graph theorem tells us that some induced subgraph is an odd hole or the complement of one, and that is all.
But what if $\chi(G)$ is much larger than $\omega(G)$?
For instance, there are graphs $G$ with $\omega(G)=2$ and with $\chi(G)$ as large as we want; what induced subgraphs 
must such a graph contain?

Let $\mathcal{F}$ be a set of graphs. We say that $\mathcal{F}$ is {\em $\chi$-bounding} if for every integer
$W$, every graph $G$ with $\omega(G)\le W$ and $\chi(G)$ sufficiently large (depending on $W$) has an induced subgraph isomorphic
to a member of $\mathcal{F}$. 
Which sets of graphs $\mathcal{F}$ have this property?

This is not an easy question; the answer is not known even when $|\mathcal{F}|=1$. Let $\mathcal{F}=\{F\}$ say.
Erd\H{o}s~\cite{erdos} showed
that for every integer $k$, there is a graph $G$ with $\omega(G)=2$, with $\chi(G)\ge k$ and with girth at least $k$ (the {\em girth}
is the length of the shortest cycle). So by taking $k>|V(F)|$, it follows that if $\{F\}$ is $\chi$-bounding then $F$ must be a forest,
and a famous conjecture (proposed independently by Gy\'arf\'as~\cite{gyarfastree} and Sumner~\cite{sumner}) asserts that this is sufficient,
that $\{F\}$ is indeed $\chi$-bounding for every forest $F$. But even this is not solved, although it has been proved for several 
special types of forests (see for instance the papers by 
Gy\'arf\'as, Szemer\'edi and Tuza~\cite{gst},
Kierstead and Penrice~\cite{kp},
Kierstead and Zhu~\cite{kierstead}, and the first author~\cite{scott}).

The same argument shows that if $\mathcal{F}$ is finite then it must contain a forest, and if the Gy\'arf\'as-Sumner conjecture
is true then that is all we can deduce; every set containing a forest is $\chi$-bounding. But what can we say about
$\chi$-bounding sets  $\mathcal{F}$ that do not contain a forest?  $\mathcal{F}$ must be infinite,
and the result of Erd\H{o}s shows that $\mathcal{F}$ must contain graphs with arbitrarily 
large girth; so an obvious question is, what if all the graphs in $\mathcal{F}$ are cycles?
Which sets of cycles are $\chi$-bounding? 

Gy\'arf\'as~\cite{gyarfas} conjectured that each of three conditions might
be sufficient:

\begin{thm}
{\bf (Conjecture)} The set of all odd cycles of length at least five is $\chi$-bounding.
\end{thm}

\begin{thm} 
{\bf (Conjecture)} For each integer $\ell$, the set of all cycles of length at least $\ell$ is $\chi$-bounding. 
\end{thm}

\begin{thm}
{\bf (Conjecture)} For each integer $\ell$, the set of all odd cycles of length at least $\ell$ is $\chi$-bounding. 
\end{thm}

The third conjecture contains the first two; but in view of the strong perfect graph theorem, 
the first (which is a rephrasing of conjecture \ref{mainconj}) seems the most
natural and has attracted the most attention, and that is what we prove in this paper.
The third conjecture remains open, although
with Maria Chudnovsky, we have recently proved the second conjecture~\cite{longholes}.
(See also~\cite{Yes-C5} for related results and~\cite{oldcycle} for earlier work.)

On a similar theme,  Bonamy, Charbit and Thomass\'e~\cite{thomasse} recently proved the following:

\begin{thm}
There exists $n$ such that every graph with no induced odd cycle of length divisible by three (and in particular, with no triangles)
has chromatic number at most $n$.
\end{thm}

What could we hope for as a more general condition? 
Let us say that a set $F\subseteq \{3,4,\dots\}$ is {\em $\chi$-bounding} if the set of all cycles with length in $F$ is $\chi$-bounding.  
It is not true that every infinite set $F$ is $\chi$-bounding; not even if we restrict to sets with upper density 1
in the positive integers.  For instance, 
consider the following sequence of graphs: let $G_1$ be the null graph, and for each $i>1$, let $G_i$ be a graph with girth 
at least $2^{|V(G_{i-1})|}$ and chromatic number at least $i$.  Letting $F$ be the set of cycle lengths that do not occur 
in any $G_i$, we see that $F$ has upper density 1, while the family $\{G_i:i\ge 1\}$ has unbounded chromatic number.

In the positive direction, perhaps it is true that every set $F$ with strictly positive {\em lower} density is $\chi$-bounding?  
As a step towards resolving this, we propose the following strengthening of the Gy\'arf\'as conjectures:

\begin{thm}
{\bf (Conjecture) }Let $F\subseteq \{3,4,\dots\}$ be an infinite set of positive integers with bounded gaps (that is, 
there is some $K>0$ such that every set of $K$ consecutive positive integers contains an element of $F$).  Then $F$ is $\chi$-bounding.
\end{thm}

We have proved this for triangle-free graphs, in a strengthened form. The following is proved in~\cite{holeseq}:

\begin{thm}\label{holeseq}
For every integer $\ell\ge 0$ there exists $n$ such that in every graph $G$ with $\omega(G)<3$ and $\chi(G)>n$ there are $\ell$ holes 
with lengths $\ell$ consecutive integers.
\end{thm}

In the other direction, it would be interesting to know whether there is a $\chi$-bounding set with lower density 0.

\section{Using a cograph}

The proof of \ref{mainthm} needs a lemma which we prove in this section, and before that we need the following.
\begin{thm}\label{step}
Let $C$ be an $\omega(G)$-clique in a graph $G$ with $C\ne \emptyset$, and let $A\subseteq V(G)\setminus C$,
such that every vertex in $C$ has a neighbour in $A$.
Then there exist $a_1,a_2\in A$,  and $c_1,c_2\in C$, all distinct, such that $a_1$ is adjacent to $c_1$ and not to $c_2$,
and $a_2$ is adjacent to $c_2$ and not to $c_1$.
\end{thm}
\Proof
Choose a vertex $a_1\in A$ with as many neighbours in $C$ as possible. Since $C\cup \{a_1\}$
is not a clique (because $C$ is an $\omega(G)$-clique), there exists $c_2\in C$ nonadjacent to $a_1$. Choose $a_2\in A$ adjacent to $c_2$. From the
choice of $a_1$, $a_2$ does not have more neighbours in $C$ than $a_1$, and so there exists $c_1\in C$ adjacent to $a_1$ and not to $a_2$.
But then $a_1,a_2,c_1,c_2$ satisfy the theorem.~\bbox

\bigskip

If $X,Y$ are disjoint subsets of $V(G)$, we say $X$ is {\em complete} to $Y$ if every vertex in $X$ is adjacent to every vertex in $Y$, and 
$X$ is {\em anticomplete} to $Y$ if every vertex in $X$ is nonadjacent to every vertex in $Y$.
A graph $H$ is a {\em cograph} if no induced subgraph is isomorphic to a three-edge path. We need the following, 
which has been discovered independently many times (see~\cite{cographs}):

\begin{thm}\label{cograph}
If $H$ is a cograph with more than one vertex, then there is a partition of $V(H)$ into two nonempty sets $A_1,A_2$, such that either $A_1$ is complete to $A_2$,
or $A_1$ is anticomplete to $A_2$.
\end{thm}

The {\em length} of a path or cycle is the number of edges in it, and its {\em parity} 
is $0,1$ depending whether its length is even or odd.
The {\em interior} of a path $P$ is the set of vertices of $P$ incident with
two edges of $P$.
If $A\subseteq V(G)$, an {\em $A$-path} means an induced path in $G$ with distinct ends both in $A$ and with interior
in $V(G)\setminus A$.
If $X\subseteq V(G)$, the subgraph of $G$ induced on $X$ is denoted by $G[X]$. We write $\chi(X)$ for $\chi(G[X])$ when there is no danger of ambiguity.
We apply \ref{cograph} to prove the following. 

\begin{thm}\label{cographuse}
Let $W> 0$, and let $G$ be a graph with $\omega(G)\le W$. Let $A,B$ be a partition of $V(G)$, such that:
\begin{itemize}
\item $A$ is stable;
\item every vertex in $B$ has a neighbour in $A$; and
\item there is a cograph $H$ with vertex set $A$, with the property that for every $A$-path $P$ of $G$, its ends are adjacent in $H$ if and only if $P$
has odd length.
\end{itemize}
Then there is a partition $X,Y$ of $B$ such that every $W$-clique in $B$ intersects both $X$ and $Y$.
\end{thm}
\Proof
We remark that the last bullet above is a very strong statement; it implies immediately that for every two
vertices $u,v\in A$, either every $A$-path of $G$ joining $u,v$ has odd length or every such $A$-path has even length.
We prove by induction on $|A|$ that if $Z$ denotes the union of all $W$-cliques in $B$, then there is a 
partition $X,Y$ of $Z$ such that every $W$-clique in $B$ intersects both $X$ and $Y$.
\\
\\
(1) {\em If $C\subseteq B$ is a $W$-clique, and $M\subseteq A$, and every vertex in $C$ has a neighbour in $M$, then no vertex in $C$ is complete to $M$.}
\\
\\
Since $M$ is stable, \ref{step} implies that there is a four-vertex induced path $m\d c\d c'\d m'$
of $G$, with $m,m'\in M$ and $c,c'\in C$. Since this is an $A$-path of odd length, $m,m'$ are adjacent in $H$.
If $v\in C$ is complete to $M$, then 
$m\d v\d m'$ is an $A$-path of even length between $m,m'$, which is impossible.
This proves (1).

\bigskip

We may assume there is a $W$-clique included in $B$, for otherwise we can take $X=Y=\emptyset$; and so $|A|\ge 2$ (since $W>0$ and every vertex in 
this $W$-clique has a neighbour in $A$, and they do not all have the same neighbour since $\omega(G)\le W$).
By \ref{cograph}, there is a partition $A_1,A_2$ of $A$ into two nonempty
sets $A_1,A_2$, such that $A_1$ is either complete or anticomplete to $A_2$ in $H$.
For $i=1,2$ let $H_i=H[A_i]$.

Suppose first that $A_1$ is complete to $A_2$ in $H$. For $i=1,2$, let
$B_i$ be the set of vertices in $B$ with a neighbour in $A_i$. Thus
$B_1\cup B_2=B$ and $B_1\cap B_2=\emptyset$. For $i = 1,2$, let
$Z_i$ be the union of all $W$-cliques included in $B_i$, and let $Z_0$ be the union of all $W$-cliques in $B$ that intersect both $B_1,B_2$.
\\
\\
(2) {$Z_0\cap Z_1, Z_0\cap Z_2=\emptyset$.}
\\
\\
For suppose that $Z_0\cap Z_1\ne \emptyset$ say. Then there is a $W$-clique $C\subseteq Z_1$ and a vertex $b_2\in B_2$ with a neighbour in $C$.
Choose $a_2\in A_2$ adjacent to $b_2$.
Let $N$ be the set of vertices in $C$ that are adjacent to $b_2$. Since $\omega(G)\le W$, it follows that $N\ne C$. 
Let $M$ be the set of vertices in $A_1$ with a neighbour in $C\setminus N$. Thus $M\ne \emptyset$.
We claim that $M$ is complete to $N$ in $G$. For let $m\in M$, and let $v\in N$. Choose $b_1\in C\setminus N$
adjacent to $m$. 
Since every $A$-path between 
$m,a_2$ has odd length, it follows that $m\d b_1\d v\d b_2\d a_2$ is not an $A$-path, and so $m$ is adjacent to $v$. This proves that
$M$ is complete to $N$.
But every vertex in $C$ has a neighbour in $M$, contrary to (1). This proves (2).

\bigskip

From the inductive hypothesis (applied to $H_1$, $A_1$ and $Z_1$), there is a partition $X_1,Y_1$ of $Z_1$ such that every $W$-clique in $Z_1$
intersects them both. Choose $X_2, Y_2$ similarly for $Z_2$. Then $X_1\cup X_2 \cup (B_1\cap Z_0)$ and $Y_1\cup Y_2\cup (B_2\cap Z_0)$
are disjoint subsets of $B$, from (2); and we claim that every $W$-clique $C$ in $B$ meets both of them. For if $C\subseteq B_1$ then
$C$ meets both $X_1$ and $Y_1$ from the choice of $X_1,Y_1$; and similarly if $C\subseteq B_2$. If $C$ intersects both $B_1,B_2$
then $C\subseteq Z_0$, and so meets both $B_1\cap Z_0$ and $B_2\cap Z_0$. This completes the inductive proof in this case.

Thus we may assume that $A_1$ is anticomplete to $A_2$ in $H$. Let $B_0$ be the set of all vertices in
$B$ with a neighbour in $A_1$ and a neighbour in $A_2$; and for $i=1,2$, let $B_i$ be the set of vertices in $B$ such that all its
neighbours in $A$ belong to $A_i$. Thus $B_0,B_1,B_2$ are pairwise disjoint and have union $B$. We claim:
\\
\\
(3) {\em Every $W$-clique in $B$ is a subset of one of $B_1, B_2$.}
\\
\\
There is no edge between $B_1$ and $B_2$, since such an edge would give a three-edge $A$-path from $A_1$ to $A_2$,
which is impossible. Thus every clique in $B$ is a subset of one of $B_0\cup B_1, B_0\cup B_2$. Suppose that $C$ is a $W$-clique contained
in $B_0\cup B_1$, with at least one vertex in $B_0$. Choose $a_2 \in A_2$ with a neighbour in $C$, and
let $N$ be the set of neighbours of $a_2$ in $C$. Thus $N\subseteq B_0$. Since $\omega(G)=|C|$, it follows that $C\setminus N\ne
\emptyset$. Let $M$ be the set of vertices in $A_1$ with a neighbour in $C\setminus N$. We claim that $M$ is complete to $N$.
For let $m\in M$ and $v\in N$. Choose $b\in C\setminus N$ adjacent to $m$. Since $m\d b\d v\d a_2$
is not an odd $A$-path (because $m,a_2$ are nonadjacent in $H$), it follows that $m,v$ are adjacent. This proves our claim that
$M$ is complete to $N$. 
But every vertex in $C$ has a neighbour in $M$, contrary to (1).
This proves (3).

\bigskip

From the inductive hypothesis (applied to $H_1$, $A_1$ and $B_1$), the union of all $W$-cliques included in $B_1$ can be partitioned into
two sets $X_1,Y_1$, such that every $W$-clique included in $B_1$ intersects them both. Choose $X_2, Y_2$ similarly for $B_2$. Then
by (3), every $W$-clique in $B$ intersects both of $X_1\cup X_2, Y_1\cup Y_2$. This completes the inductive
proof, and hence completes the proof of \ref{cographuse}.~\bbox

\section{Obtaining the cograph}

In a graph $G$, a sequence $(L_0,L_1\l L_k)$ of disjoint subsets of $V(G)$ is called a {\em levelling} in $G$ if
$|L_0|=1$, and for $1\le i\le k$, every vertex in $L_i$ has a neighbour in $L_{i-1}$, and has no neighbour in $L_h$ for $h$ with $0\le h<i-1$.

For a fixed levelling $(L_0\l L_k)$, and for
$0\le h\le j\le k$, we say that a vertex $u\in L_h$ is an {\em ancestor} of a vertex $v\in L_j$ if there is a path between $u,v$ of length $j-h$ (which therefore
has exactly one vertex in each $L_i$ for $h\le i\le j$, and has no other vertices); and if $j=h+1$ we say that $u$ is a {\em parent} of $v$.
For $0\le i\le k$, if $u,v\in L_i$ are adjacent, we say they are {\em siblings}.

If $u,v$ are vertices of a path $P$, we denote the subpath of $P$ joining them by $u\d P\d v$. If $P$ is a path between $u,v$
and $e=vw$ is an edge with $w\notin V(P)$, we denote the path with edge-set $E(P)\cup \{e\}$ by $u\d P\d v\d w$ or just $P\d v\d w$. 
If $P$ is a path between $u,v$ and $Q$ is
a path between $v,w$, and $P\cup Q$ is a path (that is, $V(P)\cap V(Q)=\{v\}$), we sometimes denote their union by $u\d P\d v\d Q\d w$, and use a similar
notation for longer strings of concatenations of paths and edges.
In this section we prove the following:

\begin{thm}\label{cographget}
Let $(L_0\l L_k)$ be a levelling in a graph $G$. Suppose that:
\begin{itemize}
\item $L_k$ is stable;
\item for every two vertices $u,v\in L_k$, all induced paths in $G$ with ends $u,v$ and with interior in 
$L_0\cup L_1\cup\cdots\cup L_{k-1}$ have the same parity; and
\item for $1\le i\le k-1$, if $u,v\in L_i$ are siblings, then they have the same set of parents (we call this the ``parent rule'').
\end{itemize}
Let $H$ be the graph with vertex set $L_{k}$, in which distinct vertices $u,v$ are adjacent if there is an induced path
of odd length between them with interior in $L_0\cup L_1\cup\cdots\cup L_{k-1}$. Then $H$ is a cograph.
\end{thm}
\Proof We proceed by induction on $|V(G)|$. Suppose there
is a four-vertex induced path in $H$, say with vertices $a_k\d b_k\d c_k\d d_k$ in order.
From the inductive hypothesis, we may assume that 
\begin{itemize}
\item $V(G)=L_0\cup L_1\cup\cdots\cup L_k$;
\item $L_k=\{a_k,b_k,c_k,d_k\}$; and
\item for $0\le i<k$, and each $v\in L_i$, there exists $u\in L_{i+1}$ such that $v$ is the only parent of $u$ (we call such a vertex $u$ a {\em dependent} of $v$).
\end{itemize}
It follows that $|L_i|\le 4$ for $0\le i\le k$.
For distinct $u,v\in L_k$, we say $(u,v)\in L_k$ is an {\em odd pair} if $u,v$ are adjacent in $H$, that is, if
every $L_k$-path between them has odd length; and $(u,v)$ is an {\em even pair} otherwise, that is, all such paths have
even length. Thus $(a_k,b_k),(b_k,c_k),(c_k,d_k)$ are odd pairs, and 
$(a_k,c_k),(b_k,d_k),(a_k,d_k)$ are even pairs. 
We will have to check the parity of several $L_k$-paths, and here is a quick rule: the parity of an $L_k$-path
is the same as the parity of the number of edges
between siblings that it contains.

We begin with the following strengthening of the parent rule:
\\
\\
(1) {\em For $0\le i\le k$, if $u,v\in L_i$ are siblings, then they have a unique parent.}
\\
\\
Suppose not, and let $p,q$ be parents of $u,v$ (each is a parent of both,
by the parent rule).
Certainly $i<k$ since $L_k$ is stable. 
Since $u,v$ each have two parents, neither of $u,v$ is a dependent of any vertex in $L_{i-1}$, and since
$p$ has a dependent, and so does $q$, it follows that $|L_i|=4$ and $L_i=\{u,v,w,x\}$ say, where $w$ is adjacent to $p$ and not to $q$,
and $x$ is adjacent to $q$ and not to $p$. Since $w,u$ do not have the same sets of parents, it follows from the parent rule that $w,u$
are not adjacent, and similarly $w,x$ are nonadjacent to $u,v$ and to each other. Now $w$ has a dependent in $L_{i+1}$, and either 
$i+1=k$ or that vertex has a dependent, and so on; so there is a path $W$ say from $w$ to $L_k$, where each vertex is a dependent
of its predecessor. Choose paths $U,V,X$ similarly for $u,v,x$. Each vertex in $W$ different from $w$ is a dependent of some other
vertex, and so has only one parent, and the same holds for $U,V,W$. Consequently $U,V,W,X$ are pairwise disjoint. There are no edges
between $U,W$, because no vertex of $U$ has a parent in $W$ and vice versa, and no vertex in $U$ has a sibling in $W$,
by the parent rule. Similarly there are no edges between any of the paths $U,V,W,X$ except for
the edge $uv$. Let the vertices of $U,V,W,X$ in $L_k$ be $u_k, v_k, w_k, x_k$; then the pair $(u_k,w_k)$ is even, because
the union of $U,W$ and the path $w\d p\d u$ is an even $L_k$-path, and similarly the pairs $(u_k,x_k), (v_k,w_k), (v_k,x_k)$
are all even, contradicting that $|E(H)|=3$. This proves (1).

\bigskip

For $i=k-1,k-2\l 0$, let $a_{i}$ be a parent of $a_{i+1}$, and define $b_i,c_i,d_i$ similarly. Thus $a_i,b_i,c_i,d_i\in L_i$ for each $i$, and and $a_0=b_0=c_0=d_0$ (since $|L_0|=1$).
Let $A$ be the path 
$$a_k\d a_{k-1}\c a_1\d a_0,$$ 
and define $B,C,D$ similarly. For each $i$ with $0\le i\le k$, let $A_i$ be the subpath $a_k\d A\d a_i$.
Since $b_1$ is adjacent to $b_0=c_0$, there
is a maximum value of $i$ with $i\le k$ such that $b_i$ has a neighbour in $C$; and so $B_i$ is disjoint from $C$. 
Since $(b_k,c_k)$ is an odd pair, and since the only possible neighbours
of $b_i$ in $C$ are $c_{i-1}, c_i, c_{i+1}$ from the definition of a levelling, it follows that $b_i$ is adjacent to $c_i$, and there are no other edges
between the paths $B_i,C_i$. Since $b_i,c_i$ are distinct it follows that $i\ge 1$; and since $L_k$ is stable it follows that $i<k$.
By (1), $b_{i-1}=c_{i-1}$.
\\
\\
(2) {\em Some vertex of $A_{i}$ has a neighbour in $B_i\cup C_i$.}
\\
\\
For suppose not. Then in particular, $A_i$ is disjoint from $B_i, C_i$. Choose an induced path $P$ between $a_{i}$ and $b_i$ with 
interior a subset of $L_0\cup\cdots\cup L_{i-1}$.
Thus the neighbour of $b_i$ in $P$ is the unique parent $b_{i-1}=c_{i-1}$ of $b_i$.
But the paths
$A_i\d a_i\d P\d b_i\d B_i$ and $A_i\d a_i\d P \d c_{i-1} \d c_i\d C_i$ have the same length, and both are $L_k$-paths, 
and one is between the odd pair $(a_k,b_k)$, and the other
between the even pair $(a_k,c_k)$, a contradiction. This proves (2). 

\bigskip

By the same argument, it follows that some
vertex of $D_i$ has a neighbour in $B_i\cup C_i$.
Choose $h\ge i$ maximum such that $a_h$ has a neighbour in $B_i\cup C_i$, and similarly choose $j\ge i$ maximum such that $d_j$ has a neighbour in $B_i\cup C_i$.
\\
\\
(3) {\em $A_h, B_i, C_i, D_j$ are pairwise disjoint, and there is no edge between $A_h, D_j$ except possibly $a_hd_j$.}
\\
\\
For we have seen that $B_i, C_i$ are disjoint, and no vertex of $A_h$ belongs to $B_i\cup C_i$ from the maximality of $h$, and the same for $D_j$.
Suppose that $A_h,D_j$ share a vertex, and choose $g\ge \max(h,j)$ maximum such that $a_g=d_g$. There is a path
from $a_g$ to $b_g$ with interior in $L_0\cup \cdots\cup L_{g-1}$, and so there is an induced path from $a_g$ to $b_k$
with interior contained in  $L_0\cup \cdots\cup L_{g-1}\cup V(B_g)$. The union of this path 
with $A_g$ and with $D_g$ are induced paths of the same length, joining an even pair and an odd pair, which is impossible. 
This proves that $A_h, D_j$ are disjoint. 

Now suppose that there is an edge between $A_h$ and $D_j$, not between $a_h$ and $d_j$; say between $a_f$ and $d_g$, where $f>h$ and $g\ge j$. Choose such an edge
with $f+g$ maximum. Since $(a_k,d_k)$ is an even pair it follows that $g=f+1$ or $g=f-1$. If $g=f-1$, then $a_{f-1}$ has no dependent,
a contradiction.
So $g=f+1$. If $g\ne j$ then $d_{g-1}$ has no dependent, again a contradiction; so
$g=j$. Hence $d_g$ has a neighbour in one of $B_i, C_i$, and there is an induced path $P$ from $d_g$ to one of $b_k,c_k$, with vertex set included in
one of $\{d_g\}\cup B_i, \{d_g\}\cup C_i$ respectively. If $g<k$, then the $L_k$-paths $D_g\cup P$ and $A_f\d a_f\d d_g\d P$ have lengths differing by two, 
and so have the same parity, yet
one joins an even pair and the other an odd pair, a contradiction. So $g=k$, and therefore $f=k-1$, and $i\le h\le k-2$.
Since $g=j=k$, $d_k$ has a neighbour in one of $B_i,C_i$. Since $L_k$ is stable, $d_k$ is not adjacent to either of $b_k,c_k$, and so it is adjacent
to at least one of $b_{k-1},c_{k-1}$. The second is impossible since $(c_k,d_k)$ is an odd pair, so $d_k$ is adjacent to  $b_{k-1}$ and to no other vertex of
$B\cup C$. There is an induced path between $a_{k-1}, c_{k-1}$ with interior in $L_0\cup\cdots\cup L_{k-2}$.
Adding the edge $c_kc_{k-1}$ to it gives an induced path $P$ between $a_{k-1},c_k$ such that
neither of $a_k,d_k$ have neighbours in its interior; but then $a_k\d a_{k-1}\d P\d c_k$ and $d_k\d a_{k-1}\d P\d c_k$
are induced
paths of the same length, one joining an odd pair and one joining an even pair, which is impossible.
This proves (3).

\bigskip
We need to figure out the possible edges between the four paths $A_h, B_i, C_i, D_j$. Except for $b_ic_i$ every such edge is incident with $a_h$ or $d_j$; 
there is symmetry between $a_h$ and $d_j$, and the
only possible neighbours of $a_h$ in $B_i, C_i, D_j$ are
$$b_{h-1}, b_h, b_{h+1},c_{h-1}, c_h, c_{h+1}, d_j.$$
Note that it is possible that $h=k$, in which case $b_{h+1}, c_{h+1}$ are not defined. 
\\
\\
(4) {\em $a_h$ has no neighbour in $C$ except possibly $c_{h-1}$.}
\\
\\
For suppose first that $h<k$ and $a_h, c_{h+1}$ are adjacent. 
Since $h\ge i$, $b_{h+1}$ is not adjacent to $c_h$; and since $c_h$ has a dependent, which cannot be any of $a_{h+1}, b_{h+1}, c_{h+1}$, it follows that $j=h+1$, and $c_h$ is the
unique parent of $d_{h+1}$. Since $(b_k,d_k)$ is an even pair, it follows that the path $D_{h+1}\d d_{h+1}\d c_h\c c_i\d b_i\d B_i$ is not induced, and
so $d_{h+1}$
has a neighbour in $B$. But $b_{h+1}$ is not a sibling of $d_{h+1}$
by the parent rule, and $b_h$ is not a parent of $d_{h+1}$ since $d_{h+1}$ is a dependent of $c_h$; and so
$h+1<k$ and $d_{h+1}$ is adjacent to $b_{h+2}$. Consequently $b_{h+1}$ has no dependent, which is
impossible. Thus it is not the case that $h<k$ and $a_h, c_{h+1}$ are adjacent; and so 
$a_h$ has no neighbour in $C$ except possibly $c_h,c_{h-1}$. Finally, if $a_h,c_h$ are adjacent,
then $A_h\d a_h\d c_h\d C_h$ is an odd $L_k$-path joining an even pair, a contradiction. This proves (4).
\\
\\
(5) {\em $a_h$ has no neighbour in $B$.}
\\
\\
Suppose that $a_h$ has a neighbour in $B$, and so $a_h$ is adjacent to one of $b_{h-1}, b_h,b_{h+1}$. 
It is not adjacent to $b_{h+1}$ since $(a_k,b_k)$ is an odd pair, 
and for the same reason if $a_h$ is adjacent to $b_{h-1}$ then it is
also adjacent to $b_h$. On the other hand, if $a_h$ is adjacent to $b_h$ then it is adjacent to $b_{h-1}$ by the parent rule. 
Thus $a_h$ is adjacent to both $b_h,b_{h-1}$.
From its definition, $h\ge i$; suppose that $h>i$.
Then the path
$A_h\d a_h\d b_{h-1}\c b_i\d c_i\d C_i$ has odd length, and therefore is not induced, since $(a_k,c_k)$ is an even pair; and so $a_h$ has a neighbour in $C$.
By (4), $a_h$ is adjacent to $c_{h-1}$; but then the siblings $a_h, b_h$ do not have the same sets of parents, which is impossible.
Thus $h=i$.

Suppose that $d_j$ has a neighbour in $C$.
Then similarly $j=i$ and $d_i$ is adjacent to $c_i,c_{i-1}$ and to no other vertex in $B\cup C$. 
Since $(a_k,d_k)$ is an even pair
it follows that $a_i,d_i$ are not adjacent. But then
$$A_i\d a_i\d b_i\d c_i\d d_i\d D_i$$ 
is an odd $L_k$-path joining an even pair, a contradiction. 

So $d_j$ has no neighbour in $C$; and since $d_j$ has a neighbour in $B_i\cup C_i$, (4) (with $a_h,d_j$ exchanged) implies
that $d_j$ is adjacent to $b_{j-1}$ and has no other
neighbour in $B\cup C$, and therefore $j>i$ (since $b_{i-1}=c_{i-1}$ and $d_j$ has no neighbour in $C$).
The path $D_j\d d_j\d b_{j-1}\c b_i$ is induced; and if $a_i,d_j$ are nonadjacent then
the union of this path with $b_i\d a_i\d A_i$ is an odd $L_k$-path joining an even pair, which is impossible.
So $a_i,d_j$ are adjacent, and since $j>i$ it follows that $j=i+1$. 
But then
$$D_{i+1}\d d_{i+1}\d a_i\d c_{i-1}\d C_{i-1}$$
is an even $L_k$-path joining an odd pair, a contradiction. This proves (5).

\bigskip
From (4) and (5), $a_h$ is adjacent to $c_{h-1}$ and to no other vertex in $B\cup C$, and similarly 
$d_j$ is adjacent to $b_{j-1}$ and to no 
other vertex in $B\cup C$. Since $b_{i-1}=c_{i-1}$ and $a_h$ has no neighbour in $B$, it follows that $h>i$,
and similarly $j>i$. If $a_h,d_j$ are nonadjacent, let $P$ be the induced path 
between $a_h$ and $d_j$ with interior in $V(B_i\cup C_i)$. Since this path uses one edge between siblings,
its union with $A_h$ and $D_j$ is an odd $L_k$-path joining an even pair, which is impossible.
So $a_h,d_j$ are adjacent, and so $j=h+1$ or $h-1$ since $(a_k,d_k)$ is an even pair, and
from the symmetry we may assume that $j=h+1$; but then $C_{h-1}\d c_{h-1}\d a_h\d d_{h+1}\d b_h\d B_h$ is an even $L_k$-path 
joining an odd pair, a contradiction.
This completes the proof of \ref{cographget}.~\bbox

By combining this and \ref{cographuse}, we deduce:

\begin{thm}\label{cographcomb}
Let $n\ge 0$ be an integer.
Let $G$ be a graph with no odd hole and let $(L_0\l L_k)$ be a levelling in $G$. Suppose that:
\begin{itemize}
\item $L_{k-1}$ is stable;
\item for $1\le i\le k-2$, if $u,v$ are adjacent vertices in $L_i$, then they have the same sets of parents; and
\item every induced subgraph of $G[L_k]$ with clique number strictly less than $\omega(G)$ is $n$-colourable.
\end{itemize}
Then $\chi(L_k)\le 2n$.
\end{thm}
\Proof
We proceed by induction on $|L_0\cup \cdots\cup L_k|$. If $G[L_k]$ is not connected, the result follows
from the inductive hypothesis applied to the levellings $(L_0\l L_{k-1}, V(C))$, for each component $C$ of $G[L_k]$ in turn. Thus we assume that
$G[L_k]$ is connected. If some vertex $v\in L_{k-1}$ has no neighbour in $L_k$, the result follows by the inductive hypothesis applied to
the levelling $(L_0,L_1\l L_{k-1}\setminus \{v\}, L_k)$. Thus we assume that every vertex in $L_{k-1}$ has a neighbour in $L_k$.
\\
\\
(1) {\em For all distinct $u,v\in L_{k-1}$, all induced paths between $u,v$ with interior in $L_0\cup \cdots\cup L_{k-2}$ have the same parity.}
\\
\\
For $u,v$ are nonadjacent since $L_{k-1}$ is stable. Since $G[L_k]$ is connected and $u,v$ both have neighbours in $L_k$, 
there is an induced path $Q$ joining $u,v$ with interior in $L_k$. For every
induced path $P$ joining $u,v$ with interior in $L_0\cup \cdots\cup L_{k-2}$, the union of $P$ and $Q$ forms a hole, which is necessarily even; and so
$P$ has the same parity as $Q$. This proves (1).

\bigskip
Let $H$ be the graph with vertex set $L_{k-1}$, in which distinct $u,v$ are adjacent if and only if some induced path joining $u,v$ with interior in
$L_0\cup \cdots\cup L_{k-2}$ has odd length. From \ref{cographget} applied to the levelling $(L_0\l L_{k-1})$, $H$ is a cograph. From \ref{cographuse} applied to 
the sets $L_{k-1}, L_k$, there is a partition $X,Y$ of $L_k$ such that $\omega(G[X]),\omega(G[Y])<\omega(G)$, and therefore $\chi(X), \chi(Y)\le n$;
and it follows that $\chi(L_k)\le 2n$. This proves \ref{cographcomb}.~\bbox

At some cost in the chromatic number, we can do without the first bullet in \ref{cographcomb}:

\begin{thm}\label{cographult}
Let $n\ge 0$ be an integer.
Let $G$ be a graph with no odd hole, and let $(L_0\l L_k)$ be a levelling in $G$. Suppose that:
\begin{itemize}
\item for $1\le i\le k-2$, if $u,v$ are adjacent vertices in $L_i$, then they have the same sets of parents; and
\item every induced subgraph of $G$ with clique number strictly less than $\omega(G)$ is $n$-colourable.
\end{itemize}
Then $\chi(L_k)\le 4n^2\omega(G)$.
\end{thm}
\Proof
We may assume that $G[L_k]$ is connected, and for $i = k-1,k-2$ and for every vertex $v\in L_i$, there exists $u\in L_{i+1}$
such that $v$ is its only parent (or else we could just delete $v$). Now we claim that $G[L_{k-2}]$ is a cograph (thanks to
Bruce Reed for this improvement); for suppose that
there is a three-edge path $a\d b\d c\d d$ that is an induced subgraph of $G[L_{k-2}]$. By the parent rule, $a,b,c,d$ have the same
parents, and in particular some vertex $z\in L_{k-3}$ is adjacent to them all. Choose $a'\in L_{k-1}$ with no neighbour in $L_{k-2}$ except
$a$, and choose $d'$ similarly for $d$. Since $a',d'$ have neighbours in $L_k$, and $G[L_k]$ is connected, there is an induced path $P$
between $a',d'$ with interior in $L_k$. But then $a\d a' \d P\d d' \d d\d z\d a$ and $a\d a' \d P\d d' \d d\d c\d b\d a$
are holes of different parity, which is impossible. This proves that $G[L_{k-2}]$ is a cograph.

In particular, $G[L_{k-2}]$ is perfect, 
so $L_{k-2}$ can be partitioned into $\omega(G)$ stable sets 
$X_1\l X_{\omega(G)}$. For $1\le j\le \omega(G)$, let $Y_j$ be the set of vertices in $L_{k-1}$ with a neighbour in $X_j$.
From \ref{cographcomb} applied to the levelling $(L_0\l L_{k-3}, X_j, Y_j)$,
it follows that $\chi(Y_j)\le 2n$. Take a partition of $Y_j$ into $2n$ stable sets $Y^1_j\l Y^{2n}_j$, and for $1\le h\le 2n$ let $Z^h_j$ be the set of vertices
in $L_k$ with a neighbour in $Y^h_j$. From \ref{cographcomb} applied to the levelling $(L_0\l L_{k-2}, Y^h_j, Z^h_j)$, it follows that
$\chi(Z^h_j)\le 2n$. Since every vertex in $L_k$ belongs to $Z^h_j$ for some choice of $h,j$, and there are $2n\omega(G)$ possible choices of $h,j$, it follows
that $\chi(L_k)\le 4n^2\omega(G)$. This proves \ref{cographult}.~\bbox

\section{The parent rule}

A major hypothesis in \ref{cographget} and \ref{cographult} was the parent rule; and in this section we show how to arrange that it holds.
In fact our proof works under the weaker hypothesis that $G$ has no odd holes of length more than $2\ell+1$, and so may be a valuable step
in proving the strongest of the three conjectures of Gy\'arf\'as mentioned earlier, although in the present paper we only use
it for $\ell=1$. Let the {\em odd hole number} of $G$ be the length of the longest induced odd cycle (or 1, if $G$ is bipartite).
First we prove the following:

\begin{thm}\label{parentrule}
Let $G$ be a graph with odd hole number at most $2\ell+1$.
Let $(L_0,L_1\l L_k)$ be a levelling in $G$. 
Then there exists a levelling $(M_0\l M_{k})$ in $G$, such that
\begin{itemize}
\item $\chi(M_k)\ge \chi(L_k)/(3\ell+3)$; and
\item for $1\le i\le k-\max(2,\ell)$, if $u,v\in M_i$ are adjacent then they have the same sets of neighbours in $M_{i-1}$.
\end{itemize}
\end{thm}

\Proof We proceed by induction on $|V(G)|$. Thus we may assume:
\begin{itemize}
\item $V(G)=L_0\cup L_1\cup \cdots\cup L_k$;
\item $G[L_k]$ is connected; and
\item for $1\le i<k$ and every vertex $u\in L_i$, there exists $v\in L_{i+1}$ such that $u$ is its only parent.
\end{itemize}

Let $s_0\in L_0$, and for $1\le i\le k$ choose $s_i\in L_i$ such that $s_{i-1}$ is its only parent. Let $S$
be the path $s_0\d s_1\c s_k$.
Let $N(S)$ be the set of vertices of $G$ not in $S$ but with a neighbour in $S$. 
If $v\in L_i\cap N(S)$, then $v$ is adjacent to one or both of $s_i,s_{i-1}$ and has no other neighbour in $S$; because
every neighbour of $v$ belongs to one of $L_{i-1},L_i,L_{i+1}$, and $v$ is not adjacent to $s_{i+1}$ since $s_i$ is the only 
parent of $s_{i+1}$.
So every vertex in $L_i\cap N(S)$ is of one of three possible types. We wish to further classify them by the value of $i$ modulo 
$\ell+1$. Let us say the {\em type} of a vertex $v\in L_i\cap N(S)$ is the pair $(\alpha,\lambda)$ where
\begin{itemize}
\item $\alpha=1,2$ or $3$ depending whether $v$ is adjacent to $s_{i-1}$ and not to $s_i$, adjacent to both $s_i$ and $s_{i-1}$, or
adjacent to $s_i$ and not to $s_{i-1}$
\item $\lambda\in \{0\l \ell\}$ is congruent to $i$ modulo $\ell+1$.
\end{itemize}
Let us fix a type $(\alpha,\lambda)=\gamma$ say. Let $V(\gamma)$ be a minimal subset of $V(G)\setminus V(S)$ such that
\begin{itemize}
\item every vertex in $N(S)$ of type $\gamma$ belongs to $V(\gamma)$; and
\item for every vertex $v\in V(G)\setminus (V(S)\cup N(S))$, if some parent of $v$ belongs to $V(\gamma)$ then $v\in V(\gamma)$.
\end{itemize}

There are $3\ell+3$ possible types $\gamma$, and every vertex in $L_k\setminus \{s_k\}$ belongs to $V(\gamma)$ for some type $\gamma$,
so either there is a type $\gamma\ne (1,1)$ such that 
$\chi(V(\gamma)\cap L_k)\ge \chi(L_k)/(3\ell+3)$, or 
$$\chi((V(\gamma)\cap L_k)\cup \{s_k\})\ge \chi(L_k)/(3\ell+3)$$ 
when $\gamma=(1,1)$. In the latter case, since $s_k$
has no neighbour in $V(\gamma)\cap L_k$, it follows that $\chi(V(\gamma)\cap L_k)\ge
\chi(L_k)/(3\ell+3)$ anyway; so we may choose a type $\gamma$ such that 
$\chi(V(\gamma)\cap L_k)\ge \chi(L_k)/(3\ell+3)$. Let $C$ be the vertex set of a component of 
$G[(V(\gamma)\cap L_k)]$
with maximum chromatic number, so $\chi(C)\ge \chi(L_k)/(3\ell+3)$. 
Let $J_k= C$, and for $i=k-1,k-2\l 1$ choose $J_{i}\subseteq V(\gamma)\cap L_{i}$ minimal such that every vertex 
in $J_{i+1}\setminus N(S)$ has a neighbour in $J_i$. 
Consequently, for every vertex $v\in J_i$, there is a path $v=p_i\d  p_{i+1}\c p_k$ such that
\begin{itemize}
\item $p_j\in J_j$ for $i\le j\le k$
\item $p_j\notin N(S)$ for $i<j\le k$
\item $p_{j-1}$ is the only parent of $p_j$ in $J_{j-1}$ for $i<j\le k$.
\end{itemize}
We call such a path a {\em pillar} for $v$. 
\\
\\
(1) {\em For $1\le i\le k-2$, if $v\in J_i$ and $v$ is nonadjacent to $s_i$, then there is an induced path $Q_v$
between $v$ and $s_i$ of length at least $2(k-i)$, with interior in $L_{i+1}\cup \cdots\cup L_k$, such that 
no vertex in $J_i$ 
different from $v$ has a neighbour in the interior of $Q_v$.}
\\
\\
Let $P_v$ be a pillar for $v$; then none of its vertices are in $N(S)$ except possibly $v$, and since both $P_v$ and the 
path $s_i\d s_{i+1}\c s_k$ end in the connected graph $G[L_k]$, there is an induced path $Q_v$ between $v,s_i$ with vertex set 
contained in the union of the vertex sets of these two 
paths and $L_k$, using all vertices of $P_v$ and $s_i\d s_{i+1}\c s_k$ except possibly their ends in $L_k$. It follows that
$Q_v$ has length at least $2(k-i)$. If $u\in J_i\setminus\{v\}$, then since $u$ has no neighbours in $L_{i+2}\cup\cdots\cup L_k$,
and $u$ is nonadjacent to the second vertex of $Q_v$ (because $v$ is its unique parent in $J_i$) and nonadjacent to $s_{i+1}$
(because $s_i$ is its only parent), it follows that $u$ has no neighbour in the interior of $Q_v$. 
This proves (1).

\bigskip

For $1\le i\le k$ and for every vertex $v\in J_i$, either $v\in N(S)$ or it has a parent in $J_{i-1}$; 
and so there is a path $v=r_i\d r_{i-1}\c r_h$
for some $h\le i$, such that $r_j\in J_j$ for $h\le j\le i$, and $r_h\in N(S)$, and $r_j\notin N(S)$ for $h+1\le j\le i$.
Since $r_h$ has a neighbour in $S$, one of
$$r_i\d  r_{i-1}\c r_h\d s_{h-1}\d  s_h\d s_{h+1}\c s_i,$$
$$r_i\d r_{i-1}\c r_h\d s_h\d s_{h+1}\c s_i$$
is an induced path (the first if $\alpha=1$ and the second if $\alpha=2$ or $3$). We choose some such path and call it $R_v$.
Note that for all $v\in J_1\cup\cdots\cup J_k$, the path $R_v$ has even length if $\alpha=1$, and odd length otherwise.
\\
\\
(2) {\em For $0\le i\le k-l$, there is no edge with one end in $J_i\cap N(S)$ and the other
in $J_i\setminus N(S)$.}
\\
\\
For suppose that $uv$ is an edge with $u\in J_i\cap N(S)$ and $v\in J_i\setminus N(S)$. 
Since $u\in N(S)$, the path $R_u$ has length one or two, and in either case adding the edge $uv$ to it gives an induced path between
$v,s_i$, with parity different from that of $R_v$. But the union of either of these paths with $Q_v$ gives a hole of length
at least $2(k-i)+2>2\ell+1$, a contradiction. This proves (2).
\\
\\
(3) {\em For $0\le i\le k-1$, if $\alpha=2$ and $u\in J_i\cap N(S)$, then $u$ has no parent in $J_{i-1}$.}
\\
\\
Because suppose $t\in J_{i-1}$ is a parent of $u$. It follows that $t\notin N(S)$ since $i,i-1$ are not congruent modulo $\ell+1$. 
Hence $R_t$ has length at least 
$2\ell+1$ (because some vertex of $R_t$ belongs to $N(S)\cap J_h$ where $h<i-1$, and so $i-h$ is a nonzero multiple of $\ell+1$,
implying $i-1-h\ge \ell$). Also $R_t$ has odd length, since $\alpha=2$. But then the union of $R_t$ with the path
$s_{i-1}\d u\d t$ is an odd hole of length at least $2\ell+3$, which is impossible. This proves (3).
\\
\\
(4) {\em For $0\le i\le \min(k-2,k-\ell)$, 
if $u,v\in J_i$ are adjacent, then they have the same sets of parents in $J_{i-1}\cup \{s_{i-1}\}$.}
\\
\\
Let $u,v\in J_i$ be adjacent. Suppose first that $u,v\in N(S)$. If $\alpha=2$ then the claim follows from (3), so we may assume that
$\alpha=1$ or $3$. In either case $u,v$ have the same sets of parents in $V(S)$; so suppose that 
there is a vertex $t\in J_{i-1}$ adjacent to $v$
and not to $u$. Choose $w\in J_i\setminus N(S)$ such that $t$ is its unique parent in $J_{i-1}$. From 
(2), $w$ is nonadjacent to $u,v$. 
Let $P_w,P_u$ be pillars for $w,u$ respectively. There is an induced path $T$ between $u,w$ with interior in $V(P_u\cup P_w)\cup C$.
If its length is even then it includes all vertices of $P_u,P_w$ except their ends in $C$, from the definition
of ``pillar'', and so has length at least $2(k-i)$; but then the union of $T$ with the path $w\d t\d v\d u$ is an odd hole of length
at least $2(k-i)+3\ge 2\ell+3$, which is impossible. So $T$ is odd. If $\alpha=1$ then $w\d t\d R_t\d s_{i-1}\d u$
is an even path joining $u,w$ of length at least $2\ell+4$; and if $\alpha=3$ then $w\d R_w\d s_i\d u$ is an even path joining $u,w$
of length at least $2\ell+4$; and in either case the union of this path with $T$ is an odd hole of length at least $2\ell+7$, 
which is impossible.

So from (2), we may assume that $u,v\notin N(S)$.
Suppose that some vertex $t\in J_{i-1}$
is adjacent to $v$ and not to $u$. The union of $R_u$ and $Q_u$ is a hole of length at least $2(k-i)+2$, and so is even; and
hence $Q_u$ and $R_u$ have the same parity, which is the same as the parity of $R_t$. 
But $v$ has no neighbour in the interior of $Q_u$, by (1);
so $u\d Q_u\d s_i\d s_{i-1}\d R_t\d t\d v\d u$ is a hole of odd length and length
at least $2(k-i)+4$, which is impossible. This proves (4).

\bigskip
If $\alpha=1$ or $2$ let $M_i=\{s_i\}\cup J_i$ for $0\le i\le k$. We claim that 
$M_0\l M_k$ satisfies the theorem.
To see this, we must check:
\begin{itemize}
\item Every vertex $v\in M_i$
has a neighbour in $M_{i-1}$. This is true if $v=s_i$ (because then it is adjacent to $s_{i-1}$),
and if $v\in N(S)$ (because then it is adjacent to $s_{i-1}$, since $\alpha=1$ or $2$), and if $v\in J_i\setminus N(S)$
(because then it has a parent in $J_{i-1}$).
\item Every vertex $v\in M_i$ has no neighbour in $M_h$ if $h<i-1$; because $M_i\subseteq L_i$ and $M_h\subseteq L_h$.
\item If $u,v\in M_i$ are adjacent and $i\le \min(k-2,k-\ell)$ then $u,v$ have the same sets of neighbours in $M_{i-1}$.
Because if $u,v\ne s_i$ then the claim follows from (4), and if say $u=s_i$ then 
$v$ is adjacent to $s_{i-1}$ (since
$\alpha=1$ or $2$), and $\alpha=2$ (since $v$ is adjacent to $s_i$), and so $u,v$ have no other neighbours in  
$M_{i-1}$ by (3). 
\end{itemize}

We assume then that $\alpha=3$. 
Let $M_0=\{s_1\}$, $M_i=\{s_{i+1}\}\cup J_i$ for $1\le i<k$, and $M_k=J_k$, and we claim that $(M_0\l M_k)$ is the desired 
levelling. Once again we must check:
\begin{itemize}
\item Every vertex $v\in M_i$
has a neighbour in $M_{i-1}$. For if $v=s_{i+1}$ then $v$ is adjacent to $s_i$; if $v\in N(S)$ then $v$ is adjacent to $s_i$
(because $\alpha=3$); and if $v\in J_i\setminus N(S)$ then it has a neighbour in $J_{i-1}$.
\item If $v\in M_i$ then it has no neighbour in $M_h$ where $h\le i-2$. Because $v\in L_i\cup L_{i+1}$, and $M_h\subseteq L_h\cup L_{h+1}$,
so this is clear unless $h=i-2$; and in that case, $v\in L_i$ is nonadjacent to $s_{i-1}$ since $\alpha=3$. 
\item If $i\le \min(k-2,k-\ell)$ and $u,v\in M_i$ are adjacent then $u,v$ have the same neighbours in $M_{i-1}$. Since $s_{i+1}$
has no neighbour in $J_i$ it follows that $u,v\ne s_{i+1}$, so 
they have the same neighbours in $J_{i-1}$ by (4), and if one is adjacent to $s_{i}$ then they both are, by (1). 
\end{itemize}
This proves \ref{parentrule}.~\bbox

\bigskip

Now we complete the proof of \ref{mainthm}.


\bigskip

\Proof
For purposes of induction it is better to prove the slightly stronger statement 
$$\chi(G)\le \frac{2^{2^{\omega(G)+2}}}{48(\omega(G)+2)}.$$
We proceed by induction on $\omega(G)$. If $\omega(G)=1$ then $\chi(G)=1$ and the result holds; 
so we may assume that $\omega(G)>1$,
and $\chi(H)\le n$ for every induced subgraph $H$ of $G$ with $\omega(H)<\omega(G)$, where 
$$n=\frac{2^{2^{(\omega+1)}}}{48(\omega+1)}$$
(writing $\omega$ for $\omega(G)$ henceforth).
We may assume that $G$ is connected; choose a vertex $s_0$, and let $L_i$ be the set of all vertices $v$ such 
that the shortest path from $s_0$ to $v$ has $i$ edges,
for all $i\ge 0$ such that some
such vertex $v$ exists. Fix $k\ge 1$. By \ref{parentrule} with $\ell=1$ there is a levelling
$(M_0\l M_{k})$ in $G$ satisfying the conclusion of \ref{parentrule}, and 
such that $\chi(M_k)\ge \chi(L_k)/6$. But by \ref{cographult}, $\chi(M_k)\le 4n^2\omega$. Consequently $\chi(L_k)\le 24n^2\omega$.
Since this holds for all $k$, it follows that 
$$\chi(G)\le 48 n^2\omega= 48\omega \left\{ \frac{2^{2^{(\omega+1)}}}{48(\omega+1)} \right\} ^2\le \frac{2^{2^{\omega+2}}}{48(\omega+2)}.$$
This proves \ref{mainthm}.~\bbox

\bigskip

One of the curious features of the proof of the strong perfect graph theorem~\cite{CRST} was that while it proved
that every graph $G$ 
either has an odd hole or has an odd hole in its complement or admits an $\omega(G)$-colouring, there was apparently
no way to convert the proof
to a polynomial-time algorithm that actually finds one of these three things
(although a polynomial-time algorithm was found shortly
afterwards, using a different method~\cite{perfectalg}). That is not the case with the result of this paper: it can easily
be converted to a polynomial-time algorithm, which, with input a graph $G$, outputs either
\begin{itemize}
\item an odd hole in $G$, or
\item a clique $C$ in $G$ and a $2^{2^{|C|+2}}$-colouring of $G$.
\end{itemize}
We leave the details to the reader. In fact, with some work we were able to bring its running time down to $O(|V(G)|^3)\log(|V(G)|))$.

\end{document}